\documentclass[12pt,a4paper]{article}
\usepackage{amssymb,theorem}
\newtheorem{thm}{Theorem}[section]
\newtheorem{lem}[thm]{Lemma}

\newtheorem{cor}[thm]{Corollary}
\newtheorem{defn}[thm]{Definition}
{\theorembodyfont{\rmfamily}  }

\begin{document}

\title{\Large{\textbf{Actions of dense subgroups of compact groups 
and $\textrm{II}_1$-factors
 with the Haagerup property}}}
\author{Paul Jolissaint}
\maketitle
\begin{abstract}
Let $M$ be a finite von Neumann algebra with the Haagerup property, and let
$G$ be a compact group that acts continuously on $M$ and that preserves
some finite trace $\tau$. We prove that if $\Gamma$ is a countable
subgroup of $G$ which has the Haagerup property, then the crossed product
algebra $M\rtimes\Gamma$ has also the Haagerup property. In
particular, we study some ergodic, non-weakly mixing actions
of groups with the Haagerup property on finite, injective von Neumann
algebras, 
and we prove that 
the associated crossed products von Neumann algebras are 
$\textrm{II}_1$-factors with the Haagerup property. If moreover the actions
have Property $(\tau)$, then the latter factors are full.
\par\vspace{3mm}\noindent

\end{abstract}

\section{Introduction}
It is well known that the Haagerup property does not
behave well under semidirect product groups \cite{ccjjv} or under
crossed products von Neumann algebras \cite{jolHAP}. Consider then an
action $\alpha$ of a countable group $\Gamma$ on a finite von Neumann
algebra $M$, such that $\alpha$ preserves some finite trace on
$M$. Assume that both $\Gamma$ and $M$ have the Haagerup property. We
look for sufficient conditions on $\alpha$ or $\Gamma$ which 
ensure that the crossed product
von Neumann algebra $M\rtimes_\alpha\Gamma$ has the same property.
\par\vspace{3mm}
For instance, Theorem 3.2 of \cite{jolHAP} provides such conditions:
if $\Gamma$ is the middle term of a short exact sequence 
$
1\rightarrow
H\rightarrow \Gamma\rightarrow Q\rightarrow 1,
$
if $Q$ is amenable and if $M\rtimes_{\alpha|H}H$ has the Haagerup
property, then so does $M\rtimes_\alpha\Gamma$.
\par\vspace{3mm}
We propose here another condition: 
\begin{thm}
Suppose that $\Gamma$ embeds into a
compact group $G$, and that $\alpha$ is the restriction to $\Gamma$ of
a continuous action of $G$ on $M$. Then $M\rtimes_\alpha\Gamma$ has the
Haagerup property if $\Gamma$ and $M$ do.
\end{thm}
A typical situation where the above conditions are satisfied is when
$\Gamma$ is maximally almost periodic and $G$ is a compact group
containing $\Gamma$. Then the action of $G$ on itself by left translations
gives a continuous action on $M=L^{\infty}(G)$ and thus
$L^\infty(G)\rtimes\Gamma$ has the Haagerup property if $\Gamma$
does. Furthermore, taking the closure of $\Gamma$ if necessary, we can
assume that it is dense in $G$, thus the corresponding crossed product
is a type $\textrm{II}_1$ factor. Observe that such actions are
non-weakly mixing.
\par\vspace{3mm}
We give next three families of examples of such factors which are
moreover full; the first one is inspired by Chapter 7 in \cite{Lub},
and the second one, by
\cite{bekka}:
\begin{thm}
Let $\mathbb H$ be the usual Hamiltonian quaternion algebra, 
let $G=\mathbb{H}^*/Z(\mathbb{H}^*)$ be the corresponding $\mathbb
Q$-algebraic group and let $p$ be any odd prime number. Embed
$\Gamma_p:=G(\mathbb{Z}[\frac{1}{p}])$ diagonally into $G(\mathbb
R)\times G(\mathbb{Q}_p)=SO(3)\times PGL_2(\mathbb{Q}_p)$, and
consider the corresponding action $\alpha$ of $\Gamma_p$ on the
2-dimensional sphere $S^2$. Then the crossed product algebra
$L^\infty(S^2)\rtimes_\alpha\Gamma_p$ is a full $II_1$-factor with the
Haagerup property.
\end{thm}
\begin{thm}
Let $\mathbb{G}$ be the $\mathbb{Q}$-algebraic group $SO(n,1)$ or $SU(n,1)$,
let $p$ be any prime number
and denote by $G_p$ the closure of $\Gamma=\mathbb{G}(\mathbb Z)$ in 
$\mathbb{G}(\mathbb{Z}_p)$. Then $L^\infty(G_p)\rtimes\Gamma$ is a full
$II_1$-factor with the Haagerup property.
\end{thm} 
Our last class of examples uses diagonal actions of $\Gamma=SL(2,\mathbb Z)$
on products of quotient groups by some principal congruence subgroups $\Gamma(m)$:
\begin{thm}
Let $\mathfrak{m}=(m_i)_{i\geq 1}$ be a sequence of integers $2\leq
m_1<m_2<\ldots$ which are pairwise coprime. Let
$G_i=\Gamma/\Gamma(m_i)$ and let $G(\mathfrak m)=\prod_{i\geq 1}G_i$
be the associated compact group on which $\Gamma$ acts
diagonally. Then $L^{\infty}(G(\mathfrak m))\rtimes\Gamma$ is a full
$II_1$-factor with the Haagerup property.
\end{thm} 
In all these theorems, fullness comes from the strong ergodicity of the actions, and, as we
will see, it is implied by the fact $\Gamma$ has Property $(\tau)$
with respect to suitable families of subgroups. See for instance
\cite{bekka} and \cite{Lub}.
\par\vspace{3mm}
Finally, we give examples of crossed product factors that are not full
but for which central sequences are under control:
\begin{thm}
Set again $\Gamma=SL(2,\mathbb Z)$, let $Z_0$ denote its center and let
$\Lambda$ be the restricted direct product group $\bigoplus_{j\geq
  1}Z_0$. Set $\tilde{\Gamma}=\Gamma\times\Lambda$ and let
$Z=Z_0\times\Lambda$ denote the center of $\tilde{\Gamma}$. Choose a
sequence $\mathfrak{m}=(m_j)_{j\geq 1}$ as in Theorem 1.4 and assume
that $m_1\geq 3$. Let $\tilde{\Gamma}$ act on $G(\mathfrak m)=\prod_{i\geq 1}G_i$
as follows:
$$
(g,(z_j)_{j\geq 1})\cdot (x_j)_{j\geq 1}=(gz_jx_j)_{j\geq 1}.
$$
Set finally $N=L^\infty(G(\mathfrak{m}))\rtimes\tilde{\Gamma}$. Then 
$N$ is a type $II_1$ factor with the Haagerup property and with
Property gamma.
Moreover, every
central sequence in $N$ is equivalent to a central sequence
$(c_n)_{n\geq 1}$ contained in the abelian von Neumann subalgebra
$L(Z)$.
\end{thm}
Our article is organized as follows: the next section contains
preliminaries on the Haagerup property, and on
Property $(\tau)$ and its relationship to strong ergodicity. Section 3
is devoted to the proof of Theorem 1.1 and the last section contains
the proofs of Theorems 1.2, 1.3, 1.4 and 1.5.


\section{Preliminaries}
\subsection{Von Neumann algebras and the Haagerup property}

Throughout the present article, $M$, $N$, denote finite von Neumann
algebras with separable preduals, $A$ denotes preferably an abelian
von Neumann algebra, 
and $\tau$ denotes some finite,
faithful, normal, normalized trace on any of these. Such a state will
be called simply a \emph{trace}. 
\par\vspace{3mm}
We
denote by $N_{*}$ the predual of $N$, by
 $L^2(N,\tau)$ the standard Hilbert space associated with
$\tau$ and by $\xi_\tau\in L^2(N,\tau)$ the unit vector which
implements $\tau$, namely such that $\tau(x)=\langle x\xi_\tau,\xi_\tau\rangle$
for every $x\in N$. We also denote by $\Vert\cdot\Vert_{2,\tau}$ the
associated Hilbert norm on both $N$ and $L^2(N,\tau)$.
If the choice of $\tau$ is fixed and that there is no danger of
confusion, we simply write $L^2(N)$ and
$\Vert\cdot\Vert_{2}$. We also set $L^2(N,\tau)_0 :=\{\xi\in
L^2(N,\tau)\ ;\ \xi\perp\xi_\tau\}$.
\par\vspace{3mm}
 Let $\mathrm{Aut}(M,\tau)$ be the group of all $\tau$-preserving
automorphisms of $M$. It is a Polish group with respect to the
topology of pointwise $\Vert\cdot\Vert_{2,\tau}$-convergence: a sequence
$(\theta_n)$ converges to $\theta$ if and only if, for all $x\in M$,
one has $\Vert \theta_n(x)-\theta(x)\Vert_{2,\tau}\to 0$ as $n\to\infty$.
\par\vspace{3mm}
Let $\Gamma$ be a countable group. The \emph{group von Neumann algebra of
  $\Gamma$}  is denoted by $L(\Gamma)$ and it is the commutant of the right regular
  representation of $\Gamma$ on $\ell^2(\Gamma)$.
Assume $\Gamma$ acts on $N$ and that the action
$\alpha$ preserves some trace $\tau$. We briefly recall the definition
and a realization 
of the corresponding crossed product $M\rtimes_\alpha\Gamma$. 
We denote by $g\mapsto\lambda_g$ (respectively $g\mapsto\rho_g$) the
left (respectively right) regular representation on $\ell^2(\Gamma)$.
We also set $\lambda(g)=1\otimes\lambda_g$, which is a unitary
operator acting on $L^2(M,\tau)\otimes\ell^2(\Gamma)$. Then
$N:=M\rtimes_\alpha\Gamma$ is the von Neumann algebra generated by 
$M\cup\{\lambda(g)\ ;\ g\in\Gamma\}$, where $x\in M$ acts on 
$L^2(M,\tau)\otimes\ell^2(\Gamma)$ as follows:
$$
(x\cdot\xi)(g)=\alpha_{g^{-1}}(x)\xi(g)
$$
for all $\xi\in
L^2(M,\tau)\otimes\ell^2(\Gamma)$, so that
$\lambda(g)x\lambda(g^{-1})=\alpha_g(x)$ for all $g$ and $x\in M$.
In this realization, $N$ is a von Neumann subalgebra of
$M\bar{\otimes}B$ (where $B$ denotes the algebra of all linear,
bounded operators on $\ell^2(\Gamma)$), and more precisely, it is the
fixed point algebra under the action $\theta$ of $\Gamma$ defined by 
$\theta_g=\alpha_g\otimes\mathrm{Ad}(\rho_g)$. We still denote by $\tau$ the
extended trace on $N$, and we let $E_M$ denote the $\tau$-preserving
conditional expectation of $N$ onto $M$. Every operator $x\in N$
admits a ``Fourier expansion''
$\displaystyle{\sum_{g\in\Gamma}x(g)\lambda(g)}$ such that
$x(g)=E_M(x\lambda(g^{-1}))\in M$ for all $g$ and $\sum\Vert
x(g)\Vert_2^2=\Vert x\Vert_2^2$.
If $\alpha:G \rightarrow \mathrm{Aut}(M,\tau)$ is a continuous action of some group
$G$ on $M$, we also denote by $g\mapsto
\alpha_g$ the corresponding unitary representation of $G$ on
$L^2(M,\tau)$ given by
$$
\alpha_g(x\xi_\tau)=\alpha_g(x)\xi_\tau\quad \forall x\in M.
$$
The restriction to the invariant subspace $L^2(M,\tau)_0 :=\{\xi\in
L^2(M,\tau)\ ;\ \xi\perp\xi_\tau\}$ is denoted by $\alpha^0$.
\par\vspace{3mm}
If $\Phi:N\rightarrow N$ is a completely
positive map such that $\tau\circ\Phi\leq \tau$, then $\Phi$ is
automatically normal and it extends to a contraction operator $T_\Phi$
on $L^2(N,\tau)$ via the equality:
$$
T_\Phi(x\xi_\tau)=\Phi(x)\xi_\tau\quad \forall x\in N.
$$
We say that $\Phi$ is $L^2$-\emph{compact} if $T_\Phi$ is a compact operator.
 Following \cite{jolHAP}, we say that $N$ has the
\emph{Haagerup property} if there exists a 
trace $\tau$ on $N$ and a sequence $(\Phi_n)_{n\geq 1}$ of completely positive, normal maps
on $N$ which satisfy:
\begin{enumerate}
\item [(i)] $\tau\circ\Phi_n\leq \tau$ and $\Phi_n$ is $L^2$-compact
  for every $n$;
\item [(ii)] for every $x\in N$, $\Vert\Phi_n(x)-x\Vert_{2,\tau}\to 0$ as $n\to\infty$.
\end{enumerate}
In fact, it follows from Proposition 2.2 of \cite{jolHAP} and
from Corollary 2, p. 39 of \cite{JoM}
that if $N$ satisfies conditions (i) and (ii) above with respect to
$\tau$, then each $\Phi_n$ can be chosen so that
$\Phi_n(1)=1$, $\tau\circ\Phi_n=\tau$ and $T_{\Phi_n}$ is a
selfadjoint
operator. Moreover, Proposition 2.4 of \cite{jolHAP} implies that $N$ satisfies the same
conditions with respect to any other trace $\tau'$ on $M$.
\par\vspace{3mm}
The above property was introduced by M. Choda in \cite{cho3} where it
was proved that
when $N$ is the group von Neumann algebra $L(\Gamma)$ of some
countable group $\Gamma$, then $L(\Gamma)$ has the Haagerup property
if and only if there exists a sequence $(\varphi_n)_{n\geq 1}$ of
positive type, normalized functions on $\Gamma$ with the following two
properties:
\begin{enumerate}
\item [(i')] for every $n$, $\varphi_n$ tends to 0 at infinity of
  $\Gamma$;
\item [(ii')] for every $g\in \Gamma$, the sequence
  $(\varphi_n(g))_{n\geq 1}$ tends to 1 as $n\to\infty$.
\end{enumerate}
See \cite{ccjjv} for much more on the Haagerup property for
locally compact groups.
\par\vspace{3mm}
Finally, consider a type $\textrm{II}_1$ factor $N$. A \emph{central
  sequence} is a bounded sequence $(x_n)_{n\geq 1}\subset N$ such that, for
  every $x\in N$,  
$$
\lim_{n\to\infty}\Vert x_nx-xx_n\Vert_2=0.
$$
Two bounded sequences $(x_n)_{n\geq 1}$ and $(y_n)_{n\geq 1}$ in 
$N$ are \emph{equivalent} if 
$$
\lim_{n\to\infty}\Vert x_n-y_n\Vert_2=0.
$$
The factor $N$ is \emph{full} if every central sequence is
\emph{trivial}, i.e. if it is equivalent to the scalar sequence
$(\tau(x_n))_{n\geq 1}$, and it has \emph{Property gamma} (of Murray
and von Neumann) if it is not full.

\subsection{Strong ergodicity and Property $(\tau)$}

Recall from \cite{KSchmidt} that a measure-preserving action of
$\Gamma$ on a probability space $(X,\mu)$ is \emph{strongly ergodic}
if every sequence $(B_n)$ of Borel subsets of $X$ that satisfies
$$
\lim_{n\to\infty}\mu(B_n\bigtriangleup gB_n)=0\quad \forall g\in\Gamma
$$
is trivial, i.e. 
$$
\lim_{n\to\infty}\mu(B_n)(1-\mu(B_n))=0.
$$
It generalizes easily to actions on finite von Neumann algebras
\cite{cho2}: 
a $\tau$-preserving action $\alpha$ of $\Gamma$ on the finite von
Neumann algebra $M$ is \emph{strongly ergodic} if every operator-norm
bounded
sequence $(x_n)\subset M$ such that $\Vert\alpha_g(x_n)-x_n\Vert_2\to
0$ as $n\to \infty$ for every $g\in \Gamma$ is equivalent to the
scalar sequence $(\tau(x_n))$, in the sense of 
Subsection 2.1. 
\par\vspace{3mm}
In \cite{cho2}, a slightly, but strictly stronger property is considered: if
$\alpha$, $\Gamma$ and $(M,\tau)$ are as above, we say that $\alpha$
is \emph{s-strongly ergodic} if, for every sequence of unit vectors $(\xi_n)_{n\geq 1}$ in
  $L^2(M)$ that satisfy
  $\Vert\alpha_g(\xi_n)-\xi_n\Vert_2\to 0$ as $n\to \infty$ for all $g\in\Gamma$,
one has $\Vert\xi_n-\langle\xi_n,\xi_\tau\rangle\xi_\tau\Vert_2\to 0$
  as $n\to\infty$. Here are characterizations of s-strong ergodicity:
\begin{lem}
Let $\alpha$, $\Gamma$ and $(M,\tau)$ be as above. 
Then the following conditions on $\alpha$ are equivalent:
\begin{enumerate}
\item [(1)] $\alpha$ is
s-strongly ergodic;
\item [(2)] $\tau$ is the unique $\alpha$-invariant state on $M$;
\item [(3)]  there exists $\delta>0$ and a
finite subset $F$ of $\Gamma$ such that 
$$
(\star)\quad \delta^2\Vert x-\tau(x)\Vert_2^2\leq\sum_{g\in
  F}\Vert\alpha_g(x)-x\Vert_2^2
\quad\forall x\in M.
$$
\end{enumerate}
\end{lem}
\emph{Proof.} The equivalence between (1) and (2) is Theorem 2 of
\cite{cho4} and it is obvious that (3) implies (1). It remains to
prove that (1) implies (3).
\par\vspace{3mm}
Thus suppose that $\alpha$ is s-strongly ergodic and let $1\in F_1\subset
F_2\subset\ \ldots\  \Gamma$ be an exhaustive sequence of finite
subsets of $\Gamma$. If one could not find $\delta>0$ and $F$
satisfying $(\star)$, there would exist a sequence $(x_n)_{n\geq 1}\in
M$ such that $\Vert x_n\Vert_2=1$, $\tau(x_n)=0$ and
$$
\sum_{g\in F_n}\Vert\alpha_g(x_n)-x_n\Vert_2^2\leq \frac{1}{n^2}
$$
for all $n$. As $\bigcup_n F_n=\Gamma$, $\xi_n=x_n\xi_\tau$ satisfies
the condition of s-strong ergodicity, but $\Vert
\xi_n-\langle\xi_n,\xi_\tau\rangle\xi_\tau\Vert$ does not converge to 0, which
is a contradiction.
\hfill Q.E.D.
\par\vspace{3mm}\noindent
\textbf{Remark.} K. Schmidt
gives in 2.7 of \cite{KSchmidt} 
an example of a strongly ergodic action of the free group $\mathbb{F}_3$ that
has more than one invariant state, thus which is not s-strongly ergodic.
\par\vspace{3mm}
The use of s-strongly ergodic actions in the context of crossed products 
is explained in the next lemma
which is adapted from \cite{cho1}:
\begin{lem}
Let $M$ be a finite von Neumann algebra equipped with some trace $\tau$
and let $\alpha$ be a $\tau$-preserving, s-strongly ergodic and free 
action of a countable group $\Gamma$ on $M$. Denote by $Z$ the center
of $\Gamma$ and assume that 
$\Gamma/Z$ is not inner amenable.
Then every central sequence in 
the crossed product type $II_1$ factor $N=M\rtimes_\alpha\Gamma$ is
equivalent to a central sequence $(c_n)_{n\geq 1}$ contained in $L(Z)$
which satisfies: for every finite subset $K$ of $Z\setminus\{1\}$, one
has
$$
\lim_{n\to\infty}\sum_{z\in K}|c_n(z)|^2=0.
$$
In particular, if $Z$ is finite, then $N$ is a full $II_1$ factor.
\end{lem}
\emph{Proof.} Let $(x_n)_{n\geq 1}\subset N$
be a central sequence. One can assume that $\Vert x_n\Vert_2=1$ for
every $n$. Let $\sum_g x_n(g)\lambda(g)$ be the Fourier expansion of
$x_n$. Then we have for every fixed $g\in \Gamma$:
\begin{eqnarray*}
 \sum_{h\in \Gamma}\vert\Vert x_n(ghg^{-1})\Vert_2 -
\Vert x_n(h)\Vert_2\vert^2 & \leq &
\sum_{h\in \Gamma}\Vert
x_n(ghg^{-1})-\lambda(g)x_n(h)\lambda(g^{-1})\Vert_2^2\\
& = & \sum_{h\in \Gamma}\Vert
\alpha_{g^{-1}}(x_n(ghg^{-1}))-x_n(h)\Vert_2^2\\
& = & \Vert x_n\lambda(g)-\lambda(g)x_n\Vert_2^2\to_{n\to\infty}0. 
\end{eqnarray*}
Since $\Gamma/Z$ is not inner amenable, one has
$$
\sum_{g\notin Z}\Vert x_n(g)\Vert_2^2\to_{n\to\infty}0.
$$
This implies that $(x_n)$ is equivalent to its projection
$(E_{M\rtimes Z}(x_n))$ onto the von Neumann subalgebra
$M\rtimes Z$. One assumes then that $x_n(g)=0$ for every $g\notin
Z$. Set $Z^*=Z\setminus\{1\}$ and let $\delta>0$ and the finite set
$F\subset\Gamma$ be as in Lemma 2.1. Then
\begin{eqnarray*}
  \sum_{g\in F}\Vert\lambda(g)x_n\lambda(g^{-1})-x_n\Vert_2^2 & = &
\sum_{z\in Z}\sum_{g\in F}\Vert\alpha_g(x_n(z))-x_n(z)\Vert_2^2 \\
 & \geq & \delta^2\sum_{z\in Z}\Vert x_n(z)-\tau(x_n(z))\Vert_2^2.
\end{eqnarray*}
Thus, $(x_n)$ is equivalent to $(c_n)\subset L(Z)$ where
$c_n=E_{L(Z)}(x_n)$ and hence $c_n(z)=\tau(x_n(z))$ for all $n$ and
$z$.
\par
Fix next a non empty finite subset $K$ of $Z^*$ and set $k=|K|$. 
For each $z\in K$, $\alpha_z$ is a properly outer automorphism, hence
there exists a non zero projection $e_z\in M$ such that 
$$
(\star)\quad\tau(e_z\alpha_z(e_z))\leq \frac{1}{2}\tau(e_z).
$$
Indeed, by Theorem 1.2.1 of \cite{co}, one takes a non zero projection $e_z$ such that $\Vert
e_z\alpha_z(e_z)\Vert\leq \frac{1}{2}$, and, as $e_z\alpha_z(e_z)\leq
\frac{1}{2}e_z$, we get $(\star)$. This implies that 
$\Vert e_z-\alpha_z(e_z)\Vert_2^2=2\tau(e_z)-2\tau(e_z\alpha_z(e_z))
\geq \tau(e_z)>0$ for every $z\in K$. Set $c=\min\{\tau(e_z)\ ;\ z\in K\}>0$
and choose a finite subset $T$ of $\Gamma\setminus Z$ of cardinality $k$
such that $tZ\cap t'Z=\emptyset$ for all $t,t'\in T$, $t\not=t'$. This
is possible since $\Gamma/Z$ is infinite (being non inner
amenable). Finally, choose some bijection $z\mapsto t_z$ from $K$ onto
$T$, set $x(t_z)=e_z$ for every $z\in K$ and put $x=\sum_{t\in
  T}x(t)\lambda(t)\in N$. Observe that, for every $z\in K$, one has
$$
 \sum_{t\in T}\Vert x(t)-\alpha_z(x(t))\Vert_2^2\geq
\Vert e_z-\alpha_z(e_z)\Vert_2^2\geq c.
$$
We get then for every $n$:
\begin{eqnarray*}
  \Vert xc_n-c_nx\Vert_2^2 & = & 
\sum_{z\in Z^*}|c_n(z)|^2\sum_{t\in T}\Vert
x(t)-\alpha_z(x(t))\Vert_2^2\\
 & \geq &  c\sum_{z\in K}|c_n(z)|^2.
\end{eqnarray*}
This proves that $\sum_{z\in K}|c_n(z)|^2\to 0$ as $n\to\infty$.
\hfill Q.E.D.

\par\vspace{3mm}
We describe next how to get s-strongly ergodic actions.
Let $\Gamma$ be a countable group and let
$\mathcal{L}=(\Gamma_\iota)_{\iota\in I}$ be a family of normal subgroups
of $\Gamma$, each $\Gamma_\iota$ having finite index in $\Gamma$. 
Denote by $R(\mathcal L)$ the family of all irreducible unitary
representations $(\rho,\mathcal{H}_\rho)$ for which there exists 
$\iota\in I$
such that $\Gamma_\iota\subset \mathrm{ker}(\rho)$. In other
words, $R(\mathcal L)$ is the subset of the unitary dual
$\hat{\Gamma}$ formed by representations that factor through some
finite quotient group $\Gamma/\Gamma_\iota$. Let us recall Definition
4.3.1 of \cite{Lub}:

\begin{defn}
Let $\Gamma$ and $\mathcal L$ be as above. We say that $\Gamma$ has
Property $(\tau)$ with respect to the family $\mathcal L$ if the
trivial representation $1_\Gamma$ is isolated in $R(\mathcal L)$.
\end{defn}
This means that one can find a positive number $\varepsilon$ and a
finite subset $F$ of $\Gamma$ such that
$$
W(\varepsilon,F)\cap R(\mathcal L)=\{1_\Gamma\},
$$
where $W(\varepsilon,F)$ is the set of
$(\rho,\mathcal{H}_\rho)\in\hat{\Gamma}$ for which there exists a unit
vector $\xi\in\mathcal{H}_\rho$ such that
$$
\max_{g\in F}\Vert\rho(g)\xi-\xi\Vert\leq \varepsilon.
$$
We will use Property $(\tau)$ in two distinct situations. 
\par\vspace{3mm}
In the first one, we consider a countable subset $\mathcal{L}'=(\Gamma_i)_{i\geq 1}$ of
$\mathcal L$.
Set $X_{\mathcal {L}'}=\prod_{i\geq 1}\Gamma/\Gamma_i$ gifted with
its natural probability measure $\mu$, and with the diagonal action of $\Gamma$:
$$
g\cdot(g_i\Gamma_i)_{i\geq 1}=((gg_i)\Gamma_i)_{i\geq 1}.
$$
Set also $A=A(\mathcal
L')=L^{\infty}(X_{\mathcal{L}'},\mu)$ and denote by $\alpha$ the
corresponding action on $A(\mathcal L')$. Integration with respect to
$\mu$ defines a trace $\tau$ on $A(\mathcal L')$.
\par\vspace{3mm}
Then Property $(\tau)$ interprets in terms of the action $\alpha$ as
follows:
\begin{lem}
Let $\Gamma$, $\mathcal L$, $\mathcal{L}'$, $A(\mathcal{L}')$,
$\mu$ and
$\tau$ be as above. Assume moreover that:
\begin{enumerate}
\item [(a)] $\Gamma$ has Property $(\tau)$ with respect to $\mathcal L$;
\item [(b)] for every finite subset $\{\Gamma_{m_1},\ \ldots,\
  \Gamma_{m_n}\}$ of $\mathcal{L}'$ there exists $\iota\in I$ such
  that
$$
\Gamma_\iota\subset \bigcap_{j=1}^n \Gamma_{m_j};
$$
\item [(c)] the action of $\Gamma$ on $\prod_{i\geq
    1}\Gamma/\Gamma_{i}$ is ergodic. 
\end{enumerate}
Then $\alpha$ is s-strongly ergodic.
\end{lem}
\emph{Proof.} Since the families $\mathcal L$ and $\mathcal L'$ are fixed, we drop the
corresponding subscripts everywhere. We prove that there exist
$\delta>0$ and $F\subset \Gamma$ that satisfy $(\star)$ in Lemma 2.3.
\par\vspace{3mm}
Condition (a) implies existence of $0<\varepsilon<1/2$ and
$F$ such that $W(\varepsilon,F)\cap R=\{1_\Gamma\}$. Set
$\delta=\varepsilon/2$.
It suffices to see that $\displaystyle{\sum_{g\in F}\Vert
  \alpha_g(a)-a\Vert_2^2\geq\delta^2}$ for every $a\in A$ such that
  $\Vert a\Vert_2=1$ and $\tau(a)=0$. Suppose the contrary. There
  exists then $a\in A$, $\Vert a\Vert_2=1$ and $\tau(a)=0$ such that
$$
\sum_{g\in F}\Vert \alpha_g(a)-a\Vert_2^2<\delta^2.
$$
For $n\geq 1$, set $A_n=L^2(X_1\times\ldots\times X_n)
=L^{\infty}(X_1\times\ldots\times X_n)\subset A$,
where $X_i=\Gamma/\Gamma_i$, and let $E_n$ be the trace preserving
conditional expectation of $A$ onto $A_n$. If $n$ is large enough so
that $\Vert a-E_n(a)\Vert_2<\delta$, then we still have
$\tau(E_n(a))=0$ and $\Vert E_n(a)\Vert_2\geq 1-\delta$. Set
$$
b=\frac{E_n(a)}{\Vert E_n(a)\Vert_2}\in A_n,
$$
which satisfies $\Vert b\Vert_2=1$ and $\tau(b)=0$. Then, denoting by
$\alpha^{(n)}$ the restriction of $\alpha$ to the $\Gamma$-invariant
subalgebra $A_n$, one has:
\begin{eqnarray*}
  \sum_{g\in F}\Vert \alpha^{(n)}_g(b)-b\Vert_2^2 & = & \frac{1}{\Vert
  E_n(a)\Vert_2^2}\sum_{g\in F}\Vert E_n(\alpha_g(a)-a)\Vert_2^2\\
& \leq & \frac{\delta^2}{(1-\delta)^2}<\varepsilon^2.
\end{eqnarray*}
There exists $\iota\in I$ such that $\Gamma_\iota\subset \Gamma_i$ for
every $i\leq n$. 
As all the
$\Gamma_i$'s are normal subgroups of $\Gamma$, it follows that 
$$
\alpha^{(n)}_g(b)(x_1,\ldots,x_n)=b(x_1,\ldots,x_n)
$$
for all $(x_1,\ldots,x_n)\in\prod_{i=1}^nX_i$ and all $g\in \Gamma_\iota$.
Hence $\Gamma_\iota\subset\mathrm{ker}(\alpha^{(n)0})$, where the latter
representation is the restriction of $\alpha^{(n)}$ to
$L^2(X_1\times\ldots\times X_n)_0$. Thus, there exists an irreducible
subrepresentation $\rho$ of $\alpha^{(n)0}$ which belongs to
$W(\varepsilon,F)\cap R$. However, $\rho$ cannot be the trivial
representation since $\alpha^{(n)0}$ does not contain
$1_\Gamma$ because the action of $\Gamma$ is ergodic. This is a contradiction.
\hfill Q.E.D.
\par\vspace{3mm}
The second situation is inspired by \cite{bekka}. Suppose that $\mathcal{L}=(\Gamma_n)_{n\geq
  1}$ is a decreasing sequence of finite index normal subgroups of
  $\Gamma$ such that $\bigcap_n\Gamma_n=\{1\}$.
Let $\Gamma_c=\mathrm{proj}\lim \Gamma/\Gamma_n$ be the projective limit of
  the sequence of finite groups $(\Gamma/\Gamma_n)$ with respect to
  the natural projections $\Gamma/\Gamma_{n+1}\rightarrow
  \Gamma/\Gamma_n$. Then $\Gamma_c$ is a compact group containing $\Gamma$ as
  a dense subgroup. In fact, $\Gamma_c$ is the completion of $\Gamma$
in the topology for which the $\Gamma_n$'s form a base of
  neighbourhoods of $1$. We denote again by $\alpha$ the action of
  $\Gamma$ on $L^{\infty}(\Gamma_c)$ by left translation.
\begin{lem}
Let $\Gamma$, $\mathcal L=(\Gamma_n)$ and $\Gamma_c$ be as above. If
$\Gamma$ has Property $(\tau)$ with respect to $\mathcal L$, then the
action of $\Gamma$ on $L^{\infty}(\Gamma_c)$ is s-strongly ergodic.
\end{lem}
\emph{Proof.} For every $n$, set $A_n=L^{\infty}(\Gamma/\Gamma_n)$, so that $A_n$ is
a von Neumann subalgebra of $A$ and that $\bigcup_n A_n$ is
$\Vert\cdot \Vert_2$-dense in $A$. Let also $\alpha^{(n)}$
(respectively $\alpha^{(n)0}$) be the restriction of the action $\alpha$
to $A_n$ (respectively to $\{a\in A_n\ ;\ \tau(a)=0\}$).
\par\vspace{3mm}
Let $0<\varepsilon<\frac{1}{2}$ and $F\subset \Gamma$
finite be such that $W(\varepsilon, F)\cap R(\mathcal
L)=\{1_\Gamma\}$. Put $\delta=\varepsilon/2$. As in the proof of Lemma
2.5, assume by
contradiction that there exists $a\in A:=L^{\infty}(\Gamma_c)$ such
that $\Vert a\Vert_2=1$, $\tau(a)=0$ and
$$
\sum_{g\in F}\Vert \alpha_g(a)-a\Vert_2^2<\delta^2.
$$
By the same arguments, there exists $n$ and  $b\in A_n$ with
$\Vert b\Vert_2=1$, $\tau(b)=0$ and such that
$$
\sum_{g\in F}\Vert \alpha^{(n)}_g(b)-b\Vert_2^2
<\varepsilon^2.
$$
Hence one can find an irreducible subrepresentation $\sigma$ of
$\alpha^{(n)0}$ such that $\sigma\in W(\varepsilon, F)$. Since
$\Gamma_n =\mathrm{ker}(\alpha^{(n)})\subset \mathrm{ker}(\sigma)$, we
have $\sigma=1_\Gamma$, but this contradicts the ergodicity of
$\alpha^{(n)}$.\
\hfill\ Q.E.D.

\section{Actions of compact groups}

Let $G$ be a compact group and let $\alpha$ be a (continuous) action
of $G$ on a von Neumann algebra $M$. If the action is ergodic, it
follows from Corollary 4.2 of \cite{hkls} that $M$ is necessarily
finite and injective, and that there is a unique $G$-invariant state
on $M$ that is a trace. Even if all our examples deal with ergodic
actions, we state our main theorem for actions that are not
necessarily ergodic. 

\begin{thm}
Let $\alpha$ be a continuous action of a compact group $G$ on a
finite von Neumann algebra $M$ that preserves some trace
$\tau$. 
Assume that $G$ contains a
countable subgroup $\Gamma$ which has the Haagerup property and that 
$M$ has the same property. Then the
corresponding crossed product von Neumann algebra $N=M\rtimes_\alpha\Gamma$
has also the Haagerup property. 
\end{thm}
Proof of Theorem 2.1 follows readily from the following two lemmas:
\begin{lem}
Retain hypotheses and notations above. There exists a sequence $(\Psi_m)_{m\geq 1}$
of completely positive, unital, normal, $\tau$-preserving maps on $M$ 
with the following properties:
\begin{enumerate}
\item [(i)] $\Psi_m$ is $L^2$-compact for every $m$;
\item [(ii)] $\alpha_g\circ\Psi_m=\Psi_m\circ\alpha_g$ for all $g\in
  G$ and all $m$;
\item [(iii)] for every $x\in M$, one has $\Vert\Psi_m(x)-x\Vert_2\to
  0$ as $m\to\infty$.
\end{enumerate}
\end{lem}
\emph{Proof.} Choose a sequence $(\Phi_m)_{m\geq 1}$ of completely
positive, unital, $\tau$-preserving maps on $M$ such that the
corresponding operators $T_{\Phi_m}$ are all compact and selfadjoint,
and such that, for every $x\in M$,  $\Vert\Phi_m(x)-x\Vert_2\to 0$ as $m\to \infty$.
Define $\Psi_m$ by:
$$
\Psi_m(x)=\int_G\alpha_g\circ\Phi_m\circ\alpha_{g^{-1}}(x)dg\quad\forall
x\in M.
$$
Notice that the integral is defined in the weak sense: this means
that, for $x\in M$, $\Psi_m(x)$ is the element of $M$ characterized by 
$$
\varphi(\Psi_m(x))=\int_G\varphi(\alpha_g\circ\Phi_m\circ\alpha_{g^{-1}}(x))dg
\quad\forall \varphi\in M_{*}.
$$
Each $\Psi_m$ is a completely positive, unital, $\tau$-preserving map
on $M$, such that $\alpha_g\circ\Psi_m=\Psi_m\circ\alpha_g$ for every
$g\in G$. This proves (ii) and the first properties of the sequence
$(\Psi_m)_m$.
\par\vspace{3mm}
Let us prove that each $\Psi_m$ is $L^2$-compact. As $m$ is fixed for
the moment, put $T=T_{\Psi_m}$
and $S=T_{\Phi_m}$. Since $S$ is a selfadjoint, compact operator, there
exists a sequence $(S_k)_{k\geq 1}$ of finite-rank, selfadjoint
operators on $L^2(M)$ such that $\Vert S-S_k\Vert\to 0$,
and $\Vert S_k\Vert\leq 1$ 
$\forall k$ because $S$ itself is a contraction. 
Moreover, one has $\Vert T-\int_G\alpha_g
S_k\alpha_{g^{-1}}dg\Vert\leq \Vert S-S_k\Vert$ for every $k$. Thus,
it remains to check that 
$S_{k,G}:=\int_G\alpha_g S_k\alpha_{g^{-1}}dg$ is a compact operator. 
Let $\mathcal B$ denote the unit ball of $L^2(M)$. Since $G$ is
compact and $S_k$ is a finite-rank operator, the set 
$\Omega_k:=\{\alpha_g S_k \alpha_{g^{-1}}(\mathcal B)\ ;\ g\in G\}$ is relatively
compact. Finally, the image of $\mathcal B$ under $S_{k,G}$ is
contained in the closed convex circled hull of $\Omega_k$ which is
compact. \\
This proves claim (i).
\par\vspace{3mm}
It remains to prove
statement (iii). As the linear span of the set of projections in $M$
is norm-dense, it suffices to prove it for projections.
Thus, fix a projection $f\in M$ and $\varepsilon >0$. One has:
\begin{eqnarray*}
  \Vert\Psi_m(f)-f\Vert_2^2 & = & 
\Vert\Psi_m(f)\Vert_2^2+\Vert
f\Vert_2^2-2\mathrm{Re}\tau(\Psi_m(f)f)\\
 & \leq & 2\Vert f\Vert_2^2 -2\mathrm{Re}\tau(\Psi_m(f)f)\\
& = & 2\left(\tau(f)-\tau(\Psi_m(f)f)\right).
\end{eqnarray*}
As $G$ is a compact group, one can find a finite set $F\subset G$ such that, for
every $g\in G$, there exists $h=h(g)\in F$ that satisfies:
$$
\Vert\alpha_g(f)-\alpha_h(f)\Vert_2\leq \frac{\varepsilon^2}{6}.
$$
Furthermore, since $\Phi_m$ tends to the identity map on $M$ in the 
pointwise $\Vert\cdot\Vert_2$-topology,
there exists an integer $n$ such that
$$
\Vert \Phi_m(\alpha_h(f))-\alpha_h(f)\Vert_2\leq
\frac{\varepsilon^2}{6}\quad\forall h\in F\quad\textrm{and}\quad
\forall m\geq n.
$$
This implies that 
$$
\sup_{g\in G}\Vert \Phi_m(\alpha_g(f))\alpha_g(f)-\alpha_g(f)\Vert_2\leq
\frac{\varepsilon^2}{2} \quad \forall m\geq n.
$$
Indeed, if $m\geq n$ and if $g\in G$, let $h=h(g)\in F$ be as
above. One has:
\begin{eqnarray*}
  \Vert \Phi_m(\alpha_g(f))\alpha_g(f)-\alpha_g(f)\Vert_2 & \leq &
\Vert \Phi_m(\alpha_g(f)-\alpha_h(f))\alpha_g(f)\Vert_2 +\\
 &   &  \Vert
 \Phi_m(\alpha_h(f))\alpha_g(f)-\alpha_h(f)\alpha_g(f)\Vert_2+\\
& & \Vert \alpha_h(f)-\alpha_g(f)\Vert_2\\
& \leq & \frac{\varepsilon^2}{2}.
\end{eqnarray*}
Since 
$$
  \tau(\Psi_m(f)f)  = \int_G\tau(\Phi_m(\alpha_g(f))\alpha_g(f))dg,
$$
we get for $m\geq n$ :
\begin{eqnarray*}
  \vert\tau(f)-\tau(\Psi_m(f)f)\vert & = &
\vert \int_G\tau[\alpha_g(f)-\Phi_m(\alpha_g(f))\alpha_g(f)]dg\vert\\
& \leq & \int_G\vert
\tau[\alpha_g(f)-\Phi_m(\alpha_g(f))\alpha_g(f)]\vert dg\leq \frac{\varepsilon^2}{2}.
\end{eqnarray*}
This proves that $\Vert\Psi_m(f)-f\Vert_2\leq \varepsilon$ for all
$m\geq n$.
\hfill Q.E.D.
\par\vspace{3mm}
The next lemma is Theorem 3 in \cite{cho1}, but we sketch the proof
for the reader's convenience.

\begin{lem}
Let $M$ be a finite von Neumann algebra gifted with a normal,
faithful, normalized trace $\tau$ and let $\alpha$ be a 
$\tau$-preserving action of a countable group $\Gamma$. Assume that:
\begin{enumerate}
\item [(i)] $\Gamma$ has the Haagerup Property;
\item [(ii)] There exists a sequence $(\Psi_m)_{m\geq 1}$ of
  $\tau$-preserving, completely positive, unital, $L^2$-compact maps
  on $M$ such that $\alpha_g\circ\Psi_m=\Psi_m\circ\alpha_g$ for all
  $m$ and $g\in \Gamma$, and that $\Vert\Psi_m(x)-x\Vert_2\to 0$ as
  $m\to\infty$, for every $x\in M$. 
\end{enumerate}
Then the crossed product von
  Neumann algebra $M\rtimes\Gamma$ has the Haagerup Property.
\end{lem}
\emph{Outline of the proof.} On the one hand, let $(\varphi_n)_{n\geq
  1}$ 
be a sequence
of positive definite, normalized functions on $\Gamma$ as in the
definition of the Haagerup property for groups. Denote by $\Phi_n$ the
completely positive multiplier on the von Neumann algebra $L(\Gamma)$
associated to $\varphi_n$. 
By \cite{haa},
it extends to $M\rtimes\Gamma$ a completely
positive map still denoted by $\Phi_n$ in such a
way that 
$$
\Phi_n(x\lambda(g))=x\varphi_n(g)\lambda(g)\quad\forall x\in M,\ g\in
\Gamma.
$$
On the other hand, the restriction $\Psi'_m$ to $M\rtimes\Gamma$ of
the completely positive map $\Psi_m\otimes i_B$ on
$M\bar{\otimes}B$ has range contained in $M\rtimes\Gamma$ because 
$\theta_g\circ\Psi_m\otimes i_B=\Psi_m\otimes i_B\circ\theta_g$ for
every $g$, and $M\rtimes\Gamma = (M\bar{\otimes}B)^\theta$. It is
straightforward to check that the sequence $(\Psi'_n\circ\Phi_n)$
satisfies then all required properties to ensure that $M\rtimes\Gamma$
has the Haagerup property.
\hfill Q.E.D.
\par\vspace{3mm}\noindent
\textbf{Remark.} Let $\Gamma$ be a countable group with the Haagerup
property and let $B$ be any finite von Neumann algebra with the
Haagerup property of dimension at least
2 gifted with some trace $\tau_B$. Consider the infinite tensor
product algebra $\displaystyle{M=\bigotimes_{g\in\Gamma}(B,\tau_B)}$ on
which $\Gamma$ acts by Bernoulli shifts:
$\beta_g(\otimes_{h}x_h)=\otimes_h x_{g^{-1}h}$. $M$ has the Haagerup
property, but
we don't know whether
the corresponding crossed product $M\rtimes_\beta\Gamma$ has the
Haagerup property except if $\Gamma$ is amenable. However, M. Choda
claims on p. 88 of \cite{cho2} that $M\rtimes_\beta\Gamma$ has that
property, using Lemma 3.3, but there is a gap in her proof. Indeed,
she constructs, from completely positive, trace-preserving, unital 
$L^2$-compact maps $\Phi$ on
$B$, infinite tensor product maps $\tilde{\Phi}:=\otimes_g \Phi_g$, where
$\Phi_g=\Phi\ \forall g$. Such maps $\tilde{\Phi}$ make perfectly sense,
are completely
positive, unital, trace-preserving, but they are not $L^2$-compact in
general so Lemma 3.3 cannot be applied. 
\begin{cor}
 Let $M$ be a finite von Neumann algebra, let $\tau$ be a trace on $M$
 and let $\alpha$ be a $\tau$-preserving action of a group $\Gamma$ on
 $M$.
Assume that both $\Gamma$ and $M$ have the Haagerup property and that
 the range of $\Gamma$ through $\alpha$ is relatively compact in 
$\mathrm{Aut}(M,\tau)$. Then the crossed product
 $M\rtimes_\alpha\Gamma$ has the Haagerup property.
\end{cor}
\begin{cor}
Let $\Gamma$ be a maximally almost periodic group with the Haagerup 
property and let $G$ be a
compact group such that $\Gamma$ embeds into $G$. Then the crossed
product algebra $L^\infty(G)\rtimes\Gamma$ has the Haagerup property.
\end{cor}
\begin{cor}
Let $\Gamma$ and $G$ be as in Corollary 3.5. Assume furthermore that
$\Gamma$ is dense in $G$ and that the latter acts freely and
ergodically on a standard probability space $(X,\mu)$, and that its
action preserves $\mu$. Then the crossed product
$L^\infty(X)\rtimes\Gamma$ is a type $II_1$ factor with the Haagerup property.
\end{cor}
As is well known, the free group $\mathbb{F}_2$ embeds into $SU(2)$,
into $SO(3)$ 
and also into $SO(n+1)$ for all odd $n\geq 3$ (see \cite{DS}). These instances provide
examples of crossed products with the Haagerup property. 
In the final section we give examples of such factors that are full,
and of factors whose central sequences are under control.

\section{Examples}

Let $\Gamma$ be a countable group with the Haagerup property which is
embeddable into a compact group $G$. Thus, for our purposes, we assume that it is dense
in $G$. Consider the action $\alpha$ of $\Gamma$ by translation on
$G$. It is ergodic and free, so that the associated crossed product 
$L^\infty(G)\rtimes_\alpha \Gamma$ is a type $\textrm{II}_1$ factor.
It follows from Theorem 3.1 that it has the Haagerup property. We
present here three families of examples of pairs $(\Gamma, G)$ for which
the corresponding factor $L^\infty(G)\rtimes_\alpha \Gamma$ is a full
factor, and one family of factors that have Property gamma with some
control on central sequences.
\par\vspace{3mm}
Here is the first family of examples which is inspired by Chapter 7 of
\cite{Lub}:
Let $\mathbb{D}=\mathbb{D}(u,v)$ be a definite quaternion algebra
defined over $\mathbb Q$:
for a ring $R$, $\mathbb{D}(R)=\{x_0+x_1i+x_2j+x_3k\ ;\ x_\ell\in
R\}$ where $u$ and $v$ are rational numbers and $i^2=-u$, $j^2=-v$,
$k^2=-uv$, $ij=-ji=k$. 
(When $u,v>0$, for example, we get the standard Hamiltonian quaternion
algebra.) Let $G=\mathbb{D}^*/Z(\mathbb{D}^*)$ be the $\mathbb{Q}$-algebraic group of
invertible elements of $\mathbb{D}$ modulo the central ones, and let $p$ be
some prime number for which $\mathbb D$ splits in $\mathbb{Q}_p$ (e.g. $p$ can
be any odd prime in the case of Hamiltonian quaternions). 
Set $\Gamma_p=G(\mathbb{Z}[\frac{1}{p}])$ and embed it diagonally
into $G(\mathbb R)\times G(\mathbb{Q}_p)=SO(3)\times PGL_2(\mathbb{Q}_p)$.
Then the projection of $\Gamma_p$ to $SO(3)$ gives an embedding to a
dense subgroup of $SO(3)$. As the latter acts on the 2-sphere $S^2$,
then so does $\Gamma_p$. One has then:
\begin{thm}
With the assumptions above, the crossed product algebras
$L^\infty(SO(3))\rtimes\Gamma_p$ and $L^\infty(S^2)\rtimes\Gamma_p$
are full type $II_1$ factors with the Haagerup property.
\end{thm}
\emph{Proof.} The fact that the factors have both the Haagerup
property follows from Corollaries 3.5 and 3.6. 
Furthermore, by Section 7.2 of \cite{Lub}, the natural Lebesgue
measures on $L^\infty(SO(3))$ and on
$L^\infty(S^2)$ are the unique $\Gamma_p$-invariant means on these
algebras. In order to apply Lemmas 2.1 and 2.2, we check that
$\Gamma_p$ is not inner amenable. But the projection $\Gamma_2$ of $\Gamma_p$
into $PGL(2,\mathbb{Q}_p)$ is faithful and it is a lattice. Corollary
of Proposition 2 of \cite{HSk} implies that $\Gamma_2$, hence
$\Gamma_p$, 
is not inner amenable. In particular, its center is trivial and the
associated factors are full.
\hfill Q.E.D.
\par\vspace{3mm}
Our second family of examples is inspired by \cite{bekka}:
\begin{thm}
Let $\mathbb{G}$ be the $\mathbb{Q}$-algebraic group $SO(n,1)$ or $SU(n,1)$,
let $p$ be any prime number
and denote by $G_p$ the closure of $\Gamma=\mathbb{G}(\mathbb Z)$ in 
$\mathbb{G}(\mathbb{Z}_p)$. Then
$L^\infty(G_p)\rtimes\mathbb{G}(\mathbb Z)$ is a full
$II_1$-factor with the Haagerup property.
\end{thm}
\emph{Proof.} Set $\Gamma=\mathbb{G}(\mathbb Z)$. It 
is a lattice in a Lie group with the Haagerup
property, thus it has the same property by Theorem 4.0.1 and Proposition
6.1.5 of \cite{ccjjv}. It follows from Theorem 3.1 that
$L^\infty(G_p)\rtimes\Gamma$ is a $\textrm{II}_1$-factor with
the Haagerup property. It remains to prove that it
is a full factor. Notice first that the center of $\Gamma$ is finite and it
follows from \cite{HSk} that $\Gamma/Z(\Gamma)$ is not inner amenable.
\par
Next, for $n\geq 1$,
let $\Gamma(n)$ denote the principal congruence subgroup
$$
\Gamma(n)=\{g\in\Gamma\ ;\ g\equiv I\ (\mathrm{mod}\ n)\}.
$$
As is explained on page 509 of \cite{bekka}, $\Gamma$ has Property $(\tau)$ with
respect to the family $\mathcal{L}=(\Gamma(n))_{n\geq 1}$: this
follows from \cite{egm}, \cite{lps}, \cite{bus} when $\mathbb
G(\mathbb R)$ is
 isomorphic to $SO(n,1)$ and from \cite{li} when it is
isomorphic to $SU(n,1)$. In
particular, it also has that property with respect to the 
subfamilly
$(\Gamma(p^n))_{n\geq 1}$, and $G_p$ equals the projective limit  
$\Gamma_c=\mathrm{proj}\lim \Gamma/\Gamma(p^n)$. Lemmas 2.2
and 2.5 imply that the associated factor is full.
\hfill Q.E.D.
\par\vspace{3mm}\noindent
\textbf{Remark.} Generally, $G_p$ is different from
$\mathbb{G}(\mathbb{Z}_p)$, but it is always a finite index subgroup. 
This follows from Lemma 2 of \cite{bekka}. Notice that, 
however, $G_p=SL(2,\mathbb{Z}_p)$ in the case $\Gamma=SL(2,\mathbb Z)$.
\par\vspace{3mm}\noindent
Our last two classes of examples involve 
$\Gamma:=SL(2,\mathbb Z)$ and some of its subgroups. 
As in the proof of the above theorem, for $n\geq 1$, 
denote by $\Gamma(n)$ the
principal congruence subgroup of $\Gamma$.
Explicitely, $\Gamma(n)$ is the group of matrices 
$\displaystyle{\left(
\begin{array}{cc}
a & b \\
c & d\end{array}\right)\in \Gamma}$ such that $a\equiv d\equiv
1\pmod n$ and $b\equiv c\equiv 0\pmod n$. It is the kernel of the
natural homomorphism from $\Gamma$ onto
$SL(2,\mathbb{Z}/n\mathbb{Z})$ (see \cite{miyake},  4.2, for
instance). In particular, $\Gamma(1)=\Gamma$.
\par\vspace{3mm}
Let $\mathfrak{m}=(m_i)_{i\geq 1}$ be a sequence of integers such that
$2\leq m_1<m_2<\ldots $ and that
$(m_i,m_j)=1$ for all $i\not=j$, and let $G_j=\Gamma/\Gamma(m_j)$
be gifted with the natural action of $\Gamma$. $\Gamma$ embeds
into the compact group $G(\mathfrak{m}):=\prod_{j\geq 1}G_j$ via the
mapping
$g\mapsto (g\Gamma(m_j))_{j\geq 1}$,
because $\bigcap_j \Gamma(m_j)=\{I\}$.

\begin{lem}
  Let $\mathfrak{m}$ be as above. Then, for every integer $m_0\geq
 1$ such that $(m_0,m_j)=1$ for every $j\geq 1$, 
 the diagonal action of
  $\Gamma(m_0)$ on $G_1\times \ldots \times G_n$ is transitive for
  every $n\geq 1$. In particular, $\Gamma(m_0)$ embeds as a dense subgroup
 into  $G(\mathfrak{m})$ 
and its action 
on $L^{\infty}(G(\mathfrak{m}))$  is ergodic and free. 
\end{lem}
\emph{Proof}. Fix $n\geq 1$; it suffices to prove that the orbit of
$(\bar{1},\ldots,\bar{1})\in\prod_{j=1}^n G_j$ under the diagonal
action of $\Gamma(m_0)$ equals $\prod_{j=1}^n G_j$.
Let then $g_1,\ldots,g_n\in\Gamma$, and let us prove that there exists
$g\in\Gamma(m_0)$ such that $g\Gamma(m_j)=g_j\Gamma(m_j)$ for every
$j=1,\ldots,n$. Write 
$\displaystyle{g_j^{-1}=
\left(
\begin{array}{cc}
a_j & b_j \\
c_j & d_j\end{array}\right)}$. Then we have to find 
$\displaystyle{g=\left(
\begin{array}{cc}
x & y \\
z & t\end{array}\right)\in \Gamma(m_0)}$ such that
$$
g_j^{-1}g=
\left(
\begin{array}{cc}
a_jx+b_jz & a_jy+b_jt \\
c_jx+d_jz & c_jy+d_jt
\end{array}\right)\equiv 
\left(
\begin{array}{cc}
1 & 0 \\
0 & 1
\end{array}\right)\pmod{m_j}
$$
for every $j=0,1,\ldots,n$. Thus, we have to find integers $x,y,z,t$
such that $xt-yz=1$ and such that $x\equiv d_j$, $y\equiv -b_j$,
$z\equiv -c_j$ and $t\equiv a_j\pmod{m_j}$ for every
$j=0,1,\ldots,n$. Set $k=m_0\cdot m_1\cdots\ldots\cdot m_n$. As $(m_i,m_j)=1$ for
all $i\not=j$, it follows from the Chinese Remainder Theorem that one
can find integers $x',y',z',t'$ that are solutions mod $k$ of the
above systems, and such that $x't'-y'z'\equiv 1\pmod k$. As the
natural homomorphism $SL(2,\mathbb Z)\rightarrow
SL(2,\mathbb{Z}/k\mathbb Z)$ is onto, the existence of $x,y,z,t$ is
proved. Density of $\Gamma(m_0)$ in $G(\mathfrak{m})$ follows from the
definition of the product topology on $G(\mathfrak{m})$.
\hfill Q.E.D.
\par\vspace{3mm}
By Example 4.3.3 D in \cite{Lub}, $\Gamma=SL(2,\mathbb Z)$ has Property
$(\tau)$ with respect to the family  
$\mathcal L$ of all its 
principal congruence subgroups. Thus Lemmas 2.2, 2.3 and Theorem 3.1 give:
\begin{thm}
Let $\mathfrak{m}$ and $m_0$ be as above. Then
$L^\infty(G(\mathfrak{m}))\rtimes\Gamma(m_0)$ is a full $II_1$-factor with the
Haagerup property.
\end{thm}
\textbf{Remark.} In fact, the abelian von Neumann algebras
$L^\infty(G_p)$ in Theorem 4.2 and $L^\infty(G(\mathfrak{m}))$ in
Theorem 4.4 contain both increasing sequences of finite-dimensional,
invariant von Neumann subalgebras that are
$\Vert\cdot\Vert_2$-dense. Thus, Proposition 3.3 of \cite{jolHAP}
suffices to prove that the corresponding factors have the Haagerup property.
\par\vspace{3mm}
At last, we give a modified construction of the above framework in
order to get exemples of non full factors with controlled central
sequences. To do that, 
set again $\Gamma =SL(2,\mathbb Z)$, let $Z_0=\{I,-I\}$ denote its center and let
$\Lambda$ be the restricted direct product group $\bigoplus_{j\geq
  1}Z_0$. Set $\tilde{\Gamma}=\Gamma_0\times\Lambda$ and let
$Z=Z_0\times\Lambda$ denote the center of $\tilde{\Gamma}$. Choose a
sequence $\mathfrak{m}=(m_j)_{j\geq 1}$ as in Theorem 1.4 and assume
that $m_1\geq 3$ so that $Z\cap \Gamma(m_j)=\{I\}$ for every $j$.
Let $\tilde{\Gamma}$ act on $G(\mathfrak m)=\prod_{i\geq 1}G_i$
as follows:
$$
(g,(z_j)_{j\geq 1})\cdot (g_j\Gamma(m_j))_{j\geq 1}=(gz_jg_j\Gamma(m_j))_{j\geq 1}.
$$
It is easy to check that the action is free, and, as the action of
$\Gamma$ is s-strongly ergodic, then so is the action of $\tilde{\Gamma}$.
For future use, we define the following sequence of subsets of $Z$:
for every positive integer $k$, let $R_k$ be the set of all
$(z_j)_{j\geq 0}\in Z$ such that $z_0z_j=1$ $\forall j\leq k$; in
other words, $z=(z_j)_{j\geq 0}$ belongs to $R_k$ if and only if
either $z=(I,...,I,z_{k+1},z_{k+2},...)$ or
$z=(-I,...,-I,z_{k+1},z_{k+2},...)$.

\begin{thm}
Retain notations above and let
$N=L^\infty(G(\mathfrak{m}))\rtimes\tilde{\Gamma}$.
Then $N$ is a type $II_1$ factor with the Haagerup property and with
Property gamma. Furthermore, 
 every
central sequence in $N$ is equivalent to a central sequence
$(c_n)_{n\geq 1}$ contained in the abelian von Neumann subalgebra
$L(Z)$ which
satisfies the following condition: $(\star)$ for every  $k\geq 1$, 
$$
 \lim_{n\to\infty}\sum_{z\notin R_k}|c_n(z)|^2=0.
$$
Conversely, every bounded sequence $(c_n)_{n\geq 1}\subset L(Z)$
that satisfies $(\star)$ is a central sequence in $N$.
\end{thm}
\emph{Proof.} By Lemma 2.2, we know that every central sequence in $N$
is equivalent to a central sequence $(c_n)_{n\geq 1}$ contained in
$L(Z)$. Fix a positive integer $k$ and $x_j\in G_j$
for $j=1,...,k$. Set $a=\chi_{(x_1,...,x_k)}\in L^\infty(G_1\times
... \times G_k)\subset L^\infty(G(\mathfrak m))$. If $z\in Z$, one has
$a=\alpha_z(a)$ if and only if $z\in R_k$. Thus, 
$$
\Vert ac_n-c_na\Vert_2^2=\sum_{z\notin R_k}|c_n(z)|^2\Vert
a-\alpha_z(a)\Vert_2^2=2\sum_{z\notin R_k}|c_n(z)|^2.
$$
Hence, if $(c_n)$ is a central sequence, then it satisfies $(\star)$, and
conversely, if it satisfies $(\star)$, then it is a central sequence because
the set of all $\chi_{(x_1,...,x_k)}$ as above is total in 
$L^\infty(G(\mathfrak m))$, and each $c_n$ commutes to $\lambda(\tilde{\Gamma})$.
\hfill $\square$

\bibliographystyle{plain}
\bibliography{refactionsdiag}

\begin{thebibliography}{10}

\bibitem{bekka}
M.~Bekka.
\newblock On uniqueness of invariant means.
\newblock {\em Proc. Amer. Math. Soc.}, 126:507--514, 1998.

\bibitem{bus}
M.~Burger and P.~Sarnak.
\newblock Ramanujan duals {II}.
\newblock {\em Invent. Math.}, 106:1--11, 1993.

\bibitem{ccjjv}
P.-A. Cherix, M.~Cowling, P.~Jolissaint, P.~Julg, and A.~Valette.
\newblock {\em Groups with the {H}aagerup property ({G}romov's
  a-{T}-menability)}.
\newblock Birkh\"auser, 2001.

\bibitem{cho1}
M.~Choda.
\newblock Inner amenability and fullness.
\newblock {\em Proc. AMS}, 86:663--666, 1982.

\bibitem{cho4}
M.~Choda.
\newblock Effect of inner amenability on strong ergodicity.
\newblock {\em Math. Japonica}, 28:109--115, 1983.

\bibitem{cho3}
M.~Choda.
\newblock Group factors of the {H}aagerup type.
\newblock {\em Proc. Japan Acad.}, 59:174--209, 1983.

\bibitem{cho2}
M.~Choda.
\newblock Strong ergodicity and full $\textrm{II}_1$ factors.
\newblock In {\em Operator algebras and their connections with topology and
  ergodic theory}, volume 1132 of {\em Lecture {N}otes in {M}ath.}, pages
  84--90. Springer, {B}erlin, 1985.

\bibitem{co}
A.~Connes.
\newblock Outer conjugacy classes of automorphisms of factors.
\newblock {\em Ann. Scient. \'Ec. Norm. Sup.}, 8:383--420, 1975.

\bibitem{HSk}
P.~de~la Harpe and G.~Skandalis.
\newblock Les r\'eseaux dans les groupes semi-simples ne sont pas
  int\'erieurement moyennables.
\newblock {\em L'Ens. Math.}, 40:291--311, 1994.

\bibitem{DS}
P.~Deligne and D.~Sullivan.
\newblock Division algebras and the {H}ausdorff-{B}anach-{T}arski paradox.
\newblock {\em L'Ens. Math.}, 29:145--150, 1983.

\bibitem{egm}
J.~Elstrod, F.~Grunewald, and J.~Mennike.
\newblock Kloosterman sums for {C}lifford algebras and a lower bound for the
  positive eigenvalues of the laplacian for congruence subgroups acting on
  hyperbolic spaces.
\newblock {\em Invent. Math.}, 101:641--685, 1990.

\bibitem{haa}
U.~Haagerup.
\newblock On the dual weights for crossed products of von {N}eumann algebras
  {II}.
\newblock {\em Math. Scand.}, 43:119--140, 1978.

\bibitem{hkls}
R.~H{\o}egh-Kron, M.~B. Landstad, and E.~St{\o}rmer.
\newblock Compact ergodic groups of automorphisms.
\newblock {\em Ann. Math.}, 114:75--86, 1981.

\bibitem{jolHAP}
P.~Jolissaint.
\newblock Haagerup approximation property for finite von {N}eumann algebras.
\newblock {\em J. Operator Th.}, 48:549--571, 2002.

\bibitem{JoM}
P.~Jolissaint and F.~Martin.
\newblock Alg\`ebres de von {N}eumann finies ayant la propri\'et\'e de
  {H}aagerup et semi-groupes ${L}^2$-compacts.
\newblock {\em Bull. Belg. Math. Soc. S. Stevin}, 11:35--48, 2004.

\bibitem{li}
J.S. Li.
\newblock Kloosterman-{S}elberg zeta functions on complex hyperbolic spaces.
\newblock {\em Amer. J. Math.}, 113:653--731, 1991.

\bibitem{lps}
J.S. Li, I.I. Piatetski-Shapiro, and P.~Sarnak.
\newblock Poincar\'e series for ${SO}(n,1)$.
\newblock {\em Proc. Indian Acad. Sci. Math., Math. Sci.}, 97:231--237, 1987.

\bibitem{Lub}
A.~Lubotzky.
\newblock {\em Discrete {G}roups, {E}xpanding {G}raphs and {I}nvariant
  {M}easures}.
\newblock Birkh\"auser {V}erlag, {B}asel, 1994.

\bibitem{miyake}
T.~Miyake.
\newblock {\em Modular {F}orms}.
\newblock Springer-{V}erlag {B}erlin, 1989.

\bibitem{KSchmidt}
K.~Schmidt.
\newblock Amenability, {K}azhdan's property {T}, strong ergodicity and
  invariant means for ergodic group-actions.
\newblock {\em Ergod. Th. \& Dynam. Sys.}, 1:223--236, 1981.

\end{thebibliography}
\par
\vspace{1cm}
\noindent
\begin{flushright}
     \begin{tabular}{l}
       Institut de Math\'emathiques,\\
       Universit\'e de Neuch\^atel,\\
       Emile-Argand 11\\
       CH-2000 Neuch\^atel, Switzerland\\
       \small {paul.jolissaint@unine.ch}
     \end{tabular}
\end{flushright}

\end{document}